\documentclass[oneside, 12pt]{amsart}

\usepackage{amscd,amssymb,enumerate, amsmath,graphicx}
\usepackage{mathrsfs}
\usepackage[all]{xy}

\newcommand\mylabel[1]{\label{#1}}

\newtheorem{theorem}{Theorem}[section]
\newtheorem*{maintheorem}{Theorem}
\newtheorem{lemma}[theorem]{Lemma}
\newtheorem{proposition}[theorem]{Proposition}
\newtheorem{corollary}[theorem]{Corollary}

\theoremstyle{definition}

\theoremstyle{remark}

\DeclareFontFamily{U}{wncy}{}
\DeclareFontShape{U}{wncy}{m}{n}{<->wncyr10}{}
\DeclareSymbolFont{mcy}{U}{wncy}{m}{n}
\DeclareMathSymbol{\Sh}{\mathord}{mcy}{"58}

\newcommand{\ZZ}	{\mathbb{Z}}

\newcommand{\PP}	{\mathbb{P}}

\newcommand{\ideala}    {\mathfrak{a}}

\newcommand  {\shA}     {\mathscr{A}}

\newcommand  {\shB}     {\mathscr{B}}

\newcommand  {\shE}     {\mathscr{E}}
\newcommand  {\shF}     {\mathscr{F}}

\newcommand  {\shI}     {\mathscr{I}}

\newcommand  {\shLif}   {\mathscr{L}\!\text{\textit{if}}}

\newcommand  {\shL}     {\mathscr{L}}

\newcommand  {\foX}     {\mathfrak{X}}
\newcommand  {\foY}     {\mathfrak{Y}}

\newcommand  {\Art}     {\operatorname{Art}}

\newcommand  {\Aut}     {\operatorname{Aut}}

\newcommand  {\Div}     {\operatorname{Div}}

\newcommand  {\FEt}      {{\text{\rm FEt}}}

\newcommand  {\GL}      {\operatorname{GL}}

\newcommand  {\Hom}     {\operatorname{Hom}}

\newcommand  {\id}      {{\operatorname{id}}}

\newcommand  {\Lif}     {\operatorname{Lif}}

\newcommand  {\invlim}  {\varprojlim}

\newcommand  {\lra}     {\longrightarrow}

\newcommand  {\maxid}   {\mathfrak{m}}

\renewcommand{\O}       {\mathscr{O}}

\newcommand  {\op}      {{\operatorname{op}}}
\newcommand  {\ord}     {\operatorname{ord}}

\newcommand  {\pr}      {\operatorname{pr}}

\newcommand  {\quadand} {\quad\text{and}\quad}

\newcommand  {\ra}      {\rightarrow}

\newcommand  {\Sch}     {\text{Sch}}

\newcommand  {\Spec}    {\operatorname{Spec}}

\def\mydate{\number\day\space\ifcase\month \or January\or February\or March\or 
April\or May\or June\or July\or
August\or September\or October\or November\or December\fi \space\number\year}

\DeclareFontFamily{U}{wncy}{}
\DeclareFontShape{U}{wncy}{m}{n}{<->wncyr10}{}
\DeclareSymbolFont{mcy}{U}{wncy}{m}{n}
\DeclareMathSymbol{\Sh}{\mathord}{mcy}{"58}

\begin{document}

\title[Equivariant formal deformation theory]
      {On equivariant formal deformation theory}

\author[Stefan Schr\"oer]{Stefan Schr\"oer}
\address{Mathematisches Institut, Heinrich-Heine-Universit\"at, 40204 D\"usseldorf, Germany}
\curraddr{}
\email{schroeer@math.uni-duesseldorf.de}

\author[Yukihide Takayama]{Yukihide Takayama}     
\address{Department of Mathematical Sciences, Ritsumeikan University, 1-1-1 Nojihigashi, Kusatsu, Siga, 525-8577, Japan} 
\curraddr{} 
\email{takayama@se.ritsumei.ac.jp} 

\subjclass[2010]{14D15, 14B25, 14F17}

\dedicatory{6 April 2017}
\keywords{formal deformation theory, group actions, fibered categories, \'etale coverings, pre-Tango structures}

\begin{abstract}
Using the set-up of deformation
categories of Talpo and Vistoli, we re-interpret and generalize,  
in the context of cartesian morphisms in abstract categories, some results of Rim concerning obstructions 
against extensions of group actions in infinitesimal deformations.
Furthermore, we observe that  finite \'etale coverings can be infinitesimally extended and 
the resulting formal scheme is   algebraizable.
Finally,  we show that pre-Tango structures survive under pullbacks
with respect to finite, generically \'etale surjections $\pi:X\ra Y$,
and record some consequences regarding Kodaira vanishing in degree one.
\end{abstract}

\maketitle
\tableofcontents

\section*{Introduction}
\mylabel{Introduction}

In deformation theory, one often seeks to extend automorphisms along 
infinitesimal extensions. This is not always possible: For example,
Serre \cite{Serre 1961} showed that there are flat families of smooth hypersurfaces
$X\subset\PP^4$ over $\Lambda=\ZZ_p$ whose closed fiber $X_0$ comes with  a free
action of some elementary abelian $p$-group $G$ that does not extend
to all infinitesimal neighborhoods $X_n$. Furthermore, the resulting quotient 
$Y_0=X_0/G$ then does not lift to characteristic zero.

Rim \cite{Rim 1980} developed a formalism that explains
the \emph{obstructions} in terms of certain group cohomology in degree one and two.
Our motivation for this paper is to elucidate and perhaps simplify
Rim's arguments by extending them into a purely categorical setting,
merely using Grothendieck's notion of \emph{cartesian morphisms} 
for functors $p:\shF\ra\shE$   between arbitrary categories \cite{SGA 1},  
much in the spirit of Talpo and Vistoli \cite{Talpo; Vistoli 2013}.

Recall that a morphism  $f:\xi\ra \xi'$ in $\shF$ over a morphism $S\ra S'$ in $\shE$
is cartesian if, intuitively speaking, $\xi$ behaves like a ``base-change''    of $\xi'$ to $S$.
Now let $\xi\in \shF$ be an object over some $S\in\shE$, and  $G\ra\Aut_S(\xi)$ be a homomorphism of groups.
Write $\Lif(\xi,S')$ for the set of isomorphism classes of  cartesian morphisms
$\xi\ra\xi'$ over $S\ra S'$. This set is endowed with a $G$-action,
by transport of structure. Fix a cartesian morphism $f:\xi\ra\xi'$. Our first main result is the following:

\begin{maintheorem}
{\rm (See Theorem \ref{extension cartesian morphisms}.)}
In the above setting, suppose that the group $\Aut_\xi(\xi')$ is abelian.
Then the $G$-action on $\xi$ extends to a $G$-action on $\xi'$ if and only if
the following two conditions hold:
\begin{enumerate}
\item The isomorphism class $[f]\in\Lif(\xi,S')$ is fixed under the $G$-action.
\item The resulting cohomology class $[\tilde{G}]\in H^2(G,\Aut_\xi(\xi'))$ is trivial.
\end{enumerate}
Here $\tilde{G}=\Aut_{S'}(\xi')\times_{\Aut_S(\xi)} G$ is the induced extension of $G$
by $\Aut_\xi(\xi')$, and $[\tilde{G}]$ is the resulting cohomology class.
\end{maintheorem}

We then apply this to the following algebro-geometric setting, using the
set-up of Talpo and Vistoli \cite{Talpo; Vistoli 2013}:
Let $\Lambda$ be a complete local noetherian ring, with residue field $k=\Lambda/\maxid_\Lambda$,
and $\shF\ra (\Art_\Lambda)^\op$ be a \emph{deformation category}, that
is, a category fibered in groupoids    that satisfies 
the \emph{Rim--Schlessinger Condition}. The latter is a technical condition that comes from the structure theory of
flat schemes over Artin rings.
Note that the ring $\Lambda$ may be of mixed characteristics, which was not allowed in  \cite{Rim 1980}.
Let $\xi\in \shF(A)$ be an object, and $A'\ra A$ be a small extension of rings,
and $\xi_0=\xi|_k$.
We then use an observation of  Serre from \cite{Serre 1972} and   regard the set $\Lif(\xi,A')$, if   nonempty,
as a \emph{torsor with a group of operators} $G$, to get a a cohomology class
\begin{equation}
\label{primary obstruction}
[\Lif(\xi,A')]\in H^1(G,I\otimes_k T_{\xi_0}(\shF)).
\end{equation}
This class is trivial if and only if there is some extension $\xi\ra\xi'$ whose 
\emph{isomorphism class} is $G$-fixed. The actual $G$-action on $\xi$ extends to 
such an object  $\xi'\in\shF(A')$ if and only if the ensuing cohomology class
\begin{equation}
\label{secondary obstruction}
[\tilde{G}]\in H^2(G,\Aut_\xi(\xi'))=H^2(G,I\otimes_k\Aut_{\xi_0}(\xi_{k[\epsilon]}))
\end{equation}
vanishes. Summing up, we have a primary obstruction \eqref{primary obstruction}, 
which deals with $G$-actions on isomorphism classes,
and a secondary obstruction \eqref{secondary obstruction}, 
which takes care of the actual $G$-action on objects.

If $G$ is finite and the residue field $k=\Lambda/\maxid_\Lambda$ has characteristic $p>0$,
then the above obstructions actually lie in the corresponding
cohomology groups for a Sylow $p$-subgroup $P\subset G$.
Consequently, the $G$-action extends if and only if the $P$-action extends.

We also take up two closely related topics: First, we 
verify that finite \'etale coverings can be infinitesimally extended and 
the resulting formal scheme is always \emph{algebraizable}.
Second, we show that \emph{pre-Tango structures} survive under pullbacks
with respect to finite, generically \'etale surjections $f:X\ra Y$,
and record some consequences regarding \emph{Kodaira vanishing} in degree one.

\section{Cartesian morphisms  and extensions of group actions}
\mylabel{Cartesian morphisms}

In this section, we recall Grothendieck's notion of \emph{cartesian morphisms} (\cite{SGA 1}, Expos\'e VI),
and examine the problem of extending group actions along cartesian morphisms, using
the relation between second group cohomology and   extensions of groups.
Our motivation was to clarify and perhaps simplify some arguments of 
Rim \cite{Rim 1980}, by putting them to this  categorical  setting.

Let $p:\shF\ra\shE$ be a functor between categories $\shF$ and $\shE$.
For each object $S\in\shE$, we write $\shF(S)\subset \shF$ for the subcategory of objects $\xi$
with $p(\xi)=S$, and morphisms $h:\xi\ra\zeta$ with $p(h)=\id_S$.
The hom sets in this category are written as $\Hom_S(\xi,\zeta)$.
If $\xi\in \shF(S)$ and $\xi'\in\shF(S')$, and $S\ra S'$ is a morphism in $\shE$,
we write $\Hom_{S\ra S'}(\xi,\xi')$ for the set of morphisms $f:\xi\ra \xi'$ 
inducing the given $S\ra S'$.

Let $f:\xi\ra\xi'$ be a morphism in $\shF$, with induced morphism $S\ra S'$ in $\shE$.
One says that  $f:\xi\ra\xi'$   is   \emph{cartesian} if the map
$$
\Hom_S(\zeta,\xi) \lra \Hom_{S\ra S'}(\zeta,\xi'),\quad h\longmapsto f\circ h
$$
is bijective, for each $\zeta\in\shF(S)$. Intuitively, this means that $\xi$ is obtained from $\xi'$ by
``base-change'' along $S\ra S'$.

We also say that a cartesian morphism $f:\xi\ra\xi'$ is a \emph{lifting} of $\xi$ over $S\ra S'$.
Let $\shLif(\xi,S')$ be the set of all such liftings; by abuse of notation, we suppress the
morphism $S\ra S'$ from notation.
The group elements  $\sigma\in\Aut_S(\xi)$ act   on $\shLif(\xi,S')$ from the left
by transport of structure, written as  ${}^\sigma f = f\circ \sigma^{-1}$.
We may regard $\shLif(\xi,S')$ also as a category, where a morphism between $f:\xi\ra \xi'$
and $g:\xi\ra\zeta'$ is an $S'$-morphism $h:\xi'\ra\zeta'$ with $h\circ f= g$.
Write $\Lif(\xi,S')$ for the set of isomorphism classes $[f]$ of lifting.
Obviously, the action of $\Aut_S(\xi)$ descends to an action ${}^\sigma [f] = [f\circ\sigma^{-1}]$
from the left on $\Lif(\xi,S')$.

Every $S'$-morphism $\sigma':\xi'\ra\xi'$ yields the morphism $\sigma'\circ f$
over $S\ra S'$, which in turn corresponds to 
a unique $S$-morphism $\sigma:\xi\ra\xi$, which makes the diagram
\begin{equation}
\label{pullback morphisms}
\begin{CD}
\xi 	@>f>> 	\xi'\\
@V\sigma VV		@VV\sigma'V\\
\xi	@>>f>	\xi'
\end{CD}
\end{equation}
commutative. The map $\sigma'\mapsto\sigma$ is compatible with compositions and respects identities,
whence yields a homomorphism of groups
$$
\Aut_{S'}(\xi')\lra\Aut_S(\xi),\quad \sigma'\longmapsto\sigma.
$$
We call it the \emph{restriction map}. Its kernel $\Aut_\xi(\xi')$ equals
the group of automorphisms for the lifting $f:\xi\ra\xi'$.

Now let $G$ be a group acting on the object $\xi\in\shF(S)$, via a homomorphism of groups
$G\ra\Aut_S(\xi)$. We seek to extend this action on $\xi\in\shF(S)$ to an action on $\xi'\in\shF(S')$.
In other words, we want to complete the diagram
$$
\xymatrix{
			& G\ar[d]\ar@{-->}[dl]\\
\Aut_{S'}(\xi')\ar[r] 	& \Aut_S(\xi)
}
$$
with some dashed arrow. A necessary condition  
is that the image of $G$ in $\Aut_S(\xi)$ is contained in the image of $\Aut_{S'}(\xi')$.
This can be reformulated as a fixed point problem:

\begin{proposition}
\mylabel{fixed points}
The image of the homomorphism $G\ra\Aut_S(\xi)$ is contained in the image of $\Aut_{S'}(\xi')\ra\Aut_S(\xi)$
if and only if $[f]\in\Lif(\xi,S')$ is a fixed point for the $G$-action.
\end{proposition}

\proof
If the isomorphism class $[f]$ is fixed, then for each $\sigma\in G$, there exists an isomorphism $\sigma':\xi'\ra\xi'$
making the diagram (\ref{pullback morphisms})
commutative. Since $f$ is cartesian, the uniqueness of the arrow $\sigma$ ensures that  $\sigma'\mapsto \sigma$ under
the restriction map $\Aut_{S'}(\xi')\ra\Aut_S(\xi)$.
Conversely, if the image of $G$ lies in the image of $\Aut_{S'}(\xi')$, diagram (\ref{pullback morphisms}) shows
that the isomorphism class  of the lifting $f:\xi\ra\xi'$ is $G$-fixed.
\qed

\medskip
Now suppose that $[f]\in\Lif(\xi,S')$ is a fixed point for the $G$-action,
such that the image of $G$ in $\Aut_S(\xi)$ lies in the image of $\Aut_{S'}(\xi')$.
Setting $\tilde{G}=\Aut_{S'}(\xi')\times_{\Aut_S(\xi)} G$, we   get an induced extension of groups
\begin{equation}
\label{group extension}
1\lra\Aut_\xi(\xi')\lra \tilde{G} \lra G \lra 1.
\end{equation}
The splittings for this extension 
correspond to the extensions of the $G$-action to $\xi'$.
To express this in cohomological terms, we now make the additional assumption that the kernel 
$\Aut_\xi(\xi')$ is abelian. This abelian group becomes a $G$-module, via
${}^\sigma h = \Phi_\sigma \circ h \circ \Phi_\sigma^{-1}$,
where the $\Phi_\sigma\in\tilde{G}$ map to  $\sigma\in G$. This indeed satisfies the axioms
for actions, and does not depend on the choices of $\Phi_\sigma$, because $\Aut_\xi(\xi')$ is abelian.
Now the formula
$c_{\sigma,\tau} \Phi_{\sigma\tau} = \Phi_\sigma\Phi_\tau$ defines a cochain $c:G^2\ra\Aut_\xi(\xi')$.
As explained in \cite{Brown 1982}, Chapter IV, Section 3, this cochain is a cocycle, and the resulting cohomology class
$$
[\tilde{G}]\in H^2(G,\Aut_\xi(\xi'))
$$
does not depend on the choice  of the $\Phi_\sigma$.
Moreover, the extension of groups (\ref{group extension}) splits if and only if $[\tilde{G}]=0$.
In this case,  the extension is a semidirect product $ \Aut_\xi(\xi')\rtimes G$.
Indeed, the group $H^2(G,\Aut_\xi(\xi'))$ corresponds to isomorphism classes of group
extensions of $G$ by $\Aut_\xi(\xi')$ inducing the given $G$-module structure.
Summing up, we have shown the following ``abstract nonsense'' result:

\begin{theorem}
\mylabel{extension cartesian morphisms}
Let $p:\shF\ra\shE$ be a functor,
$f:\xi\ra\xi'$ be a cartesian morphism in $\shF$, and $S\ra S'$ be the resulting morphism in $\shE$.
Let $G\ra\Aut_S(\xi)$ be a homomorphism of groups, and assume that the group $\Aut_\xi(\xi')$ is abelian.
Then the $G$-action on $\xi\in \shF(S)$ extends to a $G$-action on $\xi'\in \shF(S')$
if and only if the following two conditions holds:
\begin{enumerate}
\item The isomorphism class $[f]\in\Lif(\xi,S')$ is fixed under the $G$-action.
\item The resulting cohomology class $[\tilde{G}]\in H^2(G,\Aut_\xi(\xi'))$ is trivial.
\end{enumerate}
\end{theorem}

Of particular practical importance  are the \emph{fibered categories} $p:\shF\ra\shE$.
This means that for each morphism $S\ra S'$ in $\shE$ and each object $\xi'\in\shF(S')$,
there is a cartesian morphism $f:\xi\ra\xi'$ in $\shF_{S\ra S'}$, and the composition of
cartesian morphisms in $\shF$ is again cartesian.
A \emph{cleavage}  is the choice, for each $S\ra S'$ and $\xi'\in\shF(S')$,
of such a cartesian morphism $f:\xi\ra\xi'$, which are called \emph{transport morphisms}.
If the transport morphisms for identities are identities, one calls the cleavage \emph{normalized}.
We   also write $ \xi'|_S=\xi$ for the domains. Intuitively, one should regard it
as a ``restriction'',   ``pull-back'' or ``base-change'' of $\xi'$ along $S\ra S'$. In fact, the transport morphisms induce
restriction or pull-back functors
$$
\shF(S')\lra\shF(S),\quad  \xi'\longmapsto \xi'|_S.
$$
In particular, for every amalgamated sum $S'\amalg_S S''$ in $\shE$, we get a functor
\begin{equation}
\label{2 product}
\shF(S'\amalg_S S'')\lra\shF(S')\times_{\shF(S)}\shF(S''), \quad \xi \longmapsto (\xi|_{S'},\xi|_{S''},\varphi),
\end{equation}
where 
$\varphi:(\xi|_{S'})|_S \lra (\xi|_{S''})|_S$
is the unique \emph{comparison isomorphism}, compare \cite{SGA 1}, Expos\'e VI, Proposition 7.2,
and the right hand side in (\ref{2 product}) is the \emph{2-fiber product of categories},
as explained in \cite{Talpo; Vistoli 2013}, Appendix C.

A \emph{category fibered in groupoids} is a fibered category $p:\shF\ra\shE$
so that the categories $\shF(S)$, with $S\in \shE$
are groupoids. 
These are the fibered categories that occur in moduli problems or deformation theory.
They have the property
that every morphism in $\shF$ is cartesian, compare \cite{SGA 1}, Expos\'e VI, Remark after Definition 6.1.

\section{Torsors with a group of operators}
\mylabel{Torsors group operators}

In this section we set up further notation, recall Serre's interpretation of first group cohomology
in terms of torsors \cite{Serre 1972}, \S 5.2, and relate it to   fixed point problems.
Let $G$ be a group that  acts from the left via automorphisms on 
another group $T$ and a set $L$. We write these actions
as $t\mapsto {}^\sigma t$ and $\xi\mapsto {}^\sigma\xi$, where $\sigma\in G$.
Suppose we have an action on the right 
$$
\mu:L\times T\lra L,\quad ( \xi,t)\longmapsto\xi\cdot t,
$$
such that the set $L$ is a \emph{principal homogeneous space} for the group $T$,
that is, a right \emph{$T$-torsor}. In other words, the set $L$ is non-empty, and for each point $\xi_0\in L$
the resulting map $T\ra L$, $t\mapsto \xi_0\cdot t$ is bijective. We assume throughout that this action is compatible with
the $G$-action in the sense
$$
{}^\sigma(\xi\cdot t ) = {}^\sigma \xi \cdot {}^\sigma t,
$$
for all $\sigma\in G$, $\xi\in L$ and $t\in T$.
One says that the $T$-torsor $L$ is \emph{endowed with a group of operators $G$}.
They are the objects of a category, where the morphisms $(L,T)\ra (L',T')$
are pairs $(f,h)$, where $f:L\ra L$ is a  $G$-equivariant maps, and $h:T\ra T'$
is a $G$-equivariant homomorphism, which satisfy
$$
f(\xi\cdot t) = f(\xi) \cdot h(t).
$$

In this situation, we want to decide whether or not the $G$-set $L$ has a fixed point.
To this end, one may construct a cohomology class $[L]\in H^1(G,T)$ as follows:
Choose some  $\xi\in L$. Then the equation ${}^\sigma \xi =   \xi\cdot t_\sigma$ defines a map
$$
G\lra T,\quad \sigma\longmapsto t_\sigma,
$$
which we regard as a $1$-cochain.  The equation
$$
\xi\cdot t_{\eta\sigma}  = 
{}^{\eta\sigma}\xi = {}^\eta (\xi\cdot t_\sigma) = 
{}^\eta \xi \cdot {}^\eta  t_\sigma =
(\xi\cdot t_\eta) \cdot {}^\eta t_\sigma  = 
\xi\cdot (t_\eta {}^\eta t_\sigma   ) 
$$
implies $t_{\eta\sigma}=t_\eta {}^\eta t_\sigma$, and it follows that the cochain is a cocycle. 
For every other point $\xi'\in L$, the equation ${}^\sigma \xi'=\xi\cdot t'_\sigma$
defines another cocycle $\sigma\mapsto t'_\sigma$.
We have $\xi'\cdot s=\xi $ for some $s \in T$, and thus 
$$
\xi'\cdot (t'_\sigma {}^\sigma s) =
{}^\sigma\xi'\cdot {}^\sigma s =
{}^\sigma (\xi'\cdot s) =
{}^\sigma\xi = 
\xi\cdot t_\sigma = 
(\xi'\cdot s ) \cdot t_\sigma =
\xi' \cdot (s  t_\sigma).
$$
It follows that $t'_\sigma = st_\sigma {}^\sigma (s^{-1})$, whence the two cocycles
are cohomologous.  We thus get a well-defined cohomology class
$$
[L]\in H^1(G,T).
$$
In this general non-abelian setting, we regard $H^1(G,T)$ as a pointed set,
where the distinguished point $\star\in H^1(G,T)$ is the cohomology class of the constant cocycle $\sigma\mapsto e$.
It is also called the \emph{trivial cohomology class}.
According to \cite{Serre 1972}, Proposition 33, this gives a pointed bijection between the set of isomorphism classes
of $T$-torsors $L$ with a group of operators $G$, and the the set $H^1(G,T)$. We need the following
consequence:

\begin{lemma}
\mylabel{cohomology fixed points}
The cohomology class $[L]\in H^1(G,T)$ is trivial if and only if the set of fixed points
$L^G$ is nonempty.
\end{lemma}

\proof
The condition is clearly sufficient: If $\xi\in L$ is $G$-fixed,
then the resulting cocycle is $t_\sigma=e$, so the cohomology class $[L]$ is trivial.
Conversely suppose that the cocycle $t_\sigma$ attached to a point $\xi\in L$ satisfies
$st_\sigma {}^\sigma (s^{-1})=e$ for some $s\in T$. Then
$$
{}^\sigma(\xi\cdot s^{-1}) = 
{}^\sigma\xi\cdot {}^\sigma s^{-1} = 
\xi \cdot t_\sigma\cdot {}^\sigma s^{-1} =
(\xi\cdot s^{-1}) \cdot (s t_\sigma{}^\sigma s^{-1}) = 
\xi\cdot s^{-1},
$$
whence $\xi'=\xi\cdot s^{-1}$ is the desired fixed point.
\qed

\section{Deformation categories and group actions}
\mylabel{Deformation categories}

Let $k$ be a field of characteristic $p\geq 0$, and let $\Lambda$ be a complete local noetherian
ring with residue field $k=\Lambda/\maxid_\Lambda$.
We write   $(\Art_\Lambda)$ for the category of
local Artin $\Lambda$-algebras $A$ such that that the induced map $k=\Lambda/\maxid_\Lambda  \ra A/\maxid_A$ on residue fields is bijective. 
Let $\shF\ra(\Art_\Lambda)^\op$
is a category fibered in groupoids satisfying the \emph{Rim--Schlessinger condition}.
Recall that the latter means that  for every cartesian square
$$
\begin{CD}
A'\times_AA'' 	@>>>	A''\\
@VVV			@VVV\\
A'		@>>>	A
\end{CD}
$$
in the category $(\Art_\Lambda)$, the resulting functor
$$
\shF(A'\times_AA'')\lra\shF(A')\times_{\shF(A)}\shF(A'')
$$
is an equivalence of categories.
Note that this functor corresponds to \eqref{2 product}, 
and is actually defined with the help  of a chosen cleavage, but the fact that it is an equivalence
does not depend on this choice.
Such a condition was first introduced by Schlessinger \cite{Schlessinger 1968}, who considered functors of Artin rings,
and extended to fibered categories by Rim \cite{Rim 1972}.
Following Talpo and Vistoli \cite{Talpo; Vistoli 2013},
we say that such a  category fibered in groupoids  $\shF\ra(\Art_\Lambda)^\op$ is a \emph{deformation category}.

Note that one should regard the opposite category $(\Art_\Lambda)^\op$ as a full subcategory
of the category $(\text{Sch}/\Lambda)$ of schemes. The morphisms in this category are thus
$\Spec(A)\ra\Spec(A')$, and correspond to algebra homomorphisms  $A'\ra A$.
The transport morphisms over a algebra homomorphism $B\ra C$, that is $\Spec(C)\ra\Spec(B)$, could also be written in
tensor product notation $\zeta\otimes_BC\ra\zeta$ instead of $\zeta|_C\ra \zeta$.
Indeed, in praxis the deformation category $\shF\ra(\Art_\Lambda)^\op$ often consists
of flat morphisms $X\ra\Spec(C)$  of  certain schemes, and the transport morphisms
are given by    projections  $\pr_1:X\otimes_BC=X\times_{\Spec(B)}\Spec(C)\ra X$.

Let $A\in (\Art_\Lambda)$, and $\xi \in \shF(A)$ be some object. Suppose that $G$ is a group
endowed with a homomorphism
$G\ra\Aut_A(\xi )$. In other words, $G$ acts on the object $\xi \in\shF$
so that the induced action on   $A \in(\Art_\Lambda)$ is trivial.
In what follows, 
$$
0\lra I\lra A'\lra A\lra 0
$$ 
is a \emph{small extension} with ideal $I\subset A'$. This means that  $I\cdot \maxid_\Lambda=0$, so we may regard the $\Lambda$-module $I$
simply as a $k$-vector space.

We now ask whether there exist a lifting $f:\xi\ra\xi'$ over $\Spec(A)\subset\Spec(A')$ to which the $G$-action
extends.
Of course, the category of all liftings may be empty, and then nothing useful can be said.
But if one assumes that some lift exists, a natural question is whether some possibly different lifting can 
be endowed with a $G$-action.
To this end, we apply Theorem \ref{extension cartesian morphisms} to our situation.
Recall that $\shLif(\xi,A')$ denotes the category of all liftings $f:\xi\ra\xi'$ over $\Spec(A)\subset\Spec(A')$,
and let $\Lif(\xi,A')$ be the set of isomorphism classes $[f]$, endowed with the canonical
$G$-action 
$$
{}^\sigma [f] = [f\circ \sigma^{-1}].
$$
To proceed, choose a  morphism $\xi_0\ra\xi$ over 
$\Spec(k)\subset\Spec(A)$, and consider the resulting \emph{tangent space}
$$
T_{\xi_0}\shF = \Lif(\xi_0,k[\epsilon]),
$$
where $\epsilon$ denotes an indeterminate subject to the relation $\epsilon^2=0$.
In other words, $k[\epsilon]\in(\Art_\Lambda)$ is the \emph{ring of dual numbers}, with ideal $k\epsilon$.

The Rim--Schlessinger condition ensures that the functor $I\mapsto \Lif(\xi_0,k[I])$ of $k$-vector spaces $I$
preserves finite products, and as a consequence $\Lif(\xi_0, k[I])$ and in particular the tangent spaces
$T_{\xi_0}\shF$ acquire the structure of an abelian group, and actually become $k$-vector spaces.
As explained in \cite{Talpo; Vistoli 2013}, Appendix A,
the natural transformation in $I$ given by
\begin{equation}
\label{tensor identification}
I\otimes_k\Lif(\xi_0,k[\epsilon])\lra \Lif(\xi_0,k[I]),\quad
v\otimes [\xi_0\stackrel{f}{\ra}\psi]\longmapsto [\xi_0\stackrel{\alpha}{\ra}\psi|_{k[I]}]
\end{equation}
is a natural isomorphism.
Here the object $ \psi|_{k[I]}$ arises from the transport morphism $\psi|_{k[I]}\ra \psi$ over the morphism
$\Spec(k[I])\ra \Spec(k[\epsilon])$
induced from the linear map 
$k\epsilon\ra I$ with $\epsilon\mapsto v$, and $\alpha:\xi_0\ra\psi|_{k[I]}$ is the transport morphism
over the inclusion $\Spec (k)\subset\Spec(k[I])$ given by $I\ra 0$.
Clearly, this natural isomorphism
respects the action of the $ \Aut_k(\xi_0)$, where the group elements $\sigma\in \Aut_k(\xi_0)$ act via
transport of structure
$$
v\otimes [\xi_0\stackrel{f}{\ra}\psi]\longmapsto v\otimes [\xi_0\stackrel{f\sigma^{-1}}{\ra}\psi]
\quadand
[\xi_0\stackrel{\alpha\sigma^{-1}}{\ra}\psi|_{k[I]}].
$$
Note that the action on $v\in I$ is trivial. In what follows, we regard the above natural isomorphism 
as an identification $I\otimes_k T_{\xi_0}(\shF)=\Lif(\xi_0,k[I])$.
Furthermore, the underlying abelian group  acts on $\Lif(\xi,A')$ in a canonical way, via some 
\begin{equation}
\label{action}
\Lif(\xi,A') \times (I\otimes_k T_{\zeta_0}\shF)\lra \Lif(\xi,A') 
\end{equation}
recalled in \eqref{new lifting} below. The $G$-action on $\xi$ induces a $G$-action on $\xi_0$,
and we also get a linear $G$-action on the tangent space $T_{\xi_0}\shF$, as described above.

\begin{proposition}
\mylabel{lifting torsor operators}
Suppose  the set  $L=\Lif(\xi,A')$ is non-empty. With respect to the   action 
of $T= I\otimes_k T_{\xi_0}\shF$, the set  $L$ is a $T$-torsor    with a group of operators $G$.
\end{proposition}

\proof
As explained in \cite{Talpo; Vistoli 2013}, Theorem 3.15,  
the Rim--Schlessinger condition ensures that the set $L$ becomes a $T$-torsor.
Our task is merely to check that this   structure is compatible with the $G$-actions. To this end, 
we have to unravel the action of $T$ on $L$.
Let $f:\xi\ra\xi'$ be lifting of $\xi\in\shF(A)$ over $A'$, and $g:\xi_0\ra \tilde{\xi}$ be a lifting
of $\xi_0\in\shF(k)$ over the ring of dual numbers $\tilde{A}=k[I]$ with ideal $I$.
We have to describe $[f]+[g]\in\Lif(\xi,A')$ and understand how the group $G$ acts on this.

To proceed, choose a cleavage for the fibered category $\shF\ra(\Art_\Lambda)^\op$. In other words,  we fix
for each object $\zeta\in\shF(C)$ and each homomorphism $B\ra C$ a transport morphism $\zeta|_C\ra \zeta$
over $\Spec(C)\ra\Spec(B)$
and regard the domain $\zeta|_C$ as the restriction of $\zeta$. We do this so that
$\xi_0=\xi|_k$ holds. In what follows, we simply write $\alpha :\zeta|_C\ra \zeta$ for 
these transport morphisms.
Now the morphism $f$ and $g$ correspond to isomorphisms
$$
\bar{f}:\xi \lra \xi'|_A\quadand \bar{g}:\xi_0\lra \tilde{\xi}|_k,
$$
and we can form the composite morphism
$$
\psi:\xi'|_k\stackrel{\bar{f}^{-1}|_k}{\lra} \xi_0\stackrel{\bar{g}}{\lra} \tilde{\xi}|_k.
$$
This gives us a triple $(\xi',\tilde{\xi},\psi)$, which we regard as an object in the fiber product category
$$
\shF(A')\times_{\shF(k)} \shF(k[I]).
$$
Now recall that we have isomorphisms of rings
$$
A'\times_A A' \lra A'\times_k (k[I]),\quad (a_1,a_2)\longmapsto (a_1,(a_1\;\text{mod}\; \maxid_{A'}, a_2-a_1)).
$$
Here we use $k[I]=k\oplus I$, and write $a_1\;\text{mod}\; \maxid_{A'}$ for the residue class in $k$,
and regard $a_2-a_1$ as element from $I$.
The Rim--Schlessinger condition yields   equivalences of categories
$$
\shF(A')\times_{\shF(A)}\shF(A')\longleftarrow \shF(A'\times_AA')\lra\shF(A'\times_kk[I]) 
\lra \shF(A')\times_{\shF(k)} \shF(k[I]),
$$
where the restriction functors are defined in terms of the chosen cleavage.
Choose adjoint equivalences, to get an equivalence of categories
$$
\shF(A')\times_{\shF(k)} \shF(k[I])\lra \shF(A')\times_{\shF(A)}\shF(A').
$$
We may choose this functor so that it commutes with the   projections onto the first factor $\shF(A')$.
Applying this functor to the object $(\xi',\tilde{\xi},\psi)$ yields an object  $(\xi',\zeta', \varphi)$,
where $\xi',\zeta'\in\shF(A')$ and $\varphi:\xi'|_A\ra\zeta'|_A$ is an isomorphism.
In turn, we get a lifting from the composite morphism  
\begin{equation}
\label{new lifting}
h:\xi\stackrel{\bar{f}}{\lra} \xi'|_A \stackrel{\varphi}{\lra} \zeta'|_A \stackrel{\alpha}{\lra} \zeta'.
\end{equation}
Here $\alpha:\zeta'|_A\ra\zeta'$ is a transport morphism.
The $T$-action on $L$ is given by $[f] + [g] = [h]$, as explained in \cite{Talpo; Vistoli 2013}, Theorem 3.15.

Now we are in the position to unravel the $G$-action.
Let $\sigma\in G$.  By definition, ${}^\sigma[f]=[f\circ\sigma^{-1}]$ and ${}^\sigma[g]=[g\circ\sigma^{-1}]$.
Using $f\circ\sigma^{-1}$ and $ g\circ\sigma^{-1}$ rather than $f$ and $g$ in the preceding paragraph, we get
$$
\overline{f\circ\sigma^{-1}}  = \bar{f}\circ\sigma^{-1}, \quadand 
\overline{g\circ\sigma^{-1}}  =   \bar{g}\circ (\sigma^{-1}|_k)=\bar{g}\circ (\sigma|_k)^{-1},
$$
which implies
$$
\overline{g\circ\sigma^{-1}} \circ \overline{f\circ\sigma^{-1}}^{-1}|_k =
\bar{g}\circ(\sigma|_k)^{-1} \circ (\sigma|_k) \circ \bar{f}^{-1}|_k = \bar{g}\circ \bar{f}^{-1}|_k.
$$
It follows that the resulting morphism $\psi:\xi'|_k\ra\tilde{\xi}$ is the same, whether computed
with $f\circ\sigma^{-1}$ and $ g\circ\sigma^{-1}$, or with $f$ and $g$.
In turn, the image of the object $(\xi',\tilde{\xi},\psi)$ remains the object $(\xi',\zeta', \varphi)$.
The resulting lifting is thus given by the composite
$$
\xi\stackrel{\bar{f}\sigma^{-1}}{\lra} \xi'|_A \stackrel{\varphi}{\lra} \zeta'|_A \stackrel{\alpha}{\lra} \zeta',
$$
which equals $h\circ\sigma^{-1}$. This shows that ${}^\sigma [f] + {}^\sigma [g] = {}^\sigma [h]$.
In other words, the $T$-torsor $L$ is endowed with a group of operators $G$.
\qed

\medskip
As described in Section \ref{Torsors group operators}, this $L$-torsor $T$ endowed with a group of operators $G$ yields
a cohomology class
$$
[\Lif(\xi,A')] \in H^1(G,I\otimes_k T_{\xi_0}\shF),
$$
and Lemma \ref{cohomology fixed points} immediately gives:

\begin{theorem}
\mylabel{cohomology fixed deformation}
Suppose   $\Lif(\xi,A')$ is non-empty.
Then the there is a $G$-fixed isomorphism class $[f]\in\Lif(\xi,A')$ 
of liftings $f:\xi \ra\xi'$  over $\Spec(A)\subset\Spec(A')$ if and only if the cohomology class
$[\Lif(\xi,A')]\in H^1(G,I\otimes_k T_{\xi_0}\shF)$ is trivial.
\end{theorem}
 
Now suppose that there exists a lifting $f:\xi\ra\xi'$ whose isomorphism class $[f]\in \Lif(\xi,A')$
is fixed under the $G$-action.
As discussed in Section \ref{Cartesian morphisms}, we get an extension of groups
\begin{equation}
\label{extension}
1\lra\Aut_{\xi}(\xi')\lra \tilde{G} \lra G\lra 1,
\end{equation}
and this extension of groups splits if and only if the $G$-action on $\xi$ extends to $\xi'$.

Now choose a morphism $\xi_{k[\epsilon]}\ra\xi_0$ over the morphism $\Spec(k[\epsilon])\ra\Spec(k)$
corresponding to the canonical inclusion $k\subset k[\epsilon]$. 
As explained in \cite{Talpo; Vistoli 2013}, Proposition 4.5,  
we have a canonical identification
\begin{equation}
\label{action automorphism}
I\otimes_k\Aut_{\xi_0}(\xi_{k[\epsilon]}) = \Aut_\xi(\xi'),
\end{equation}
and this group carries the structure of $k$-vector space. In particular, it is abelian.
In fact, \eqref{action automorphism} is    an incarnation of \eqref{action}, 
for the   deformation theory $\shA\ra(\Art_\Lambda)^\op$ whose objects over $A$ are
the automorphisms of $\xi_0|_A$, as explained in  \cite{Talpo; Vistoli 2013}, Section 4.

Since  the isomorphism class of $f:\xi\ra\xi'$ is $G$-fixed, we have a natural $G$-action
on $\Aut_{\xi'}(\xi)$, coming from the extension \eqref{extension} or equivalently from 
diagram \eqref{pullback morphisms}. The same applies for $\xi_0\ra\xi_{k[\epsilon]}$,
and we thus get a $G$-action on $\Aut_{\xi_0}(\xi_{k[\epsilon]})$.
Taking the trivial $G$-action on $I$, both sides in \eqref{action automorphism} acquire a $G$-action,
and these action coincide under the identification.
We thus may regard the extension class for (\ref{extension}) as an element in
$$
[\tilde{G}]\in H^2(G, \Aut_{\xi}(\xi')) = H^2(G,I\otimes_k\Aut_{\xi_0}(\xi_{k[\epsilon]})).
$$
Now Theorem \ref{extension cartesian morphisms} yields:

\begin{theorem}
\mylabel{cohomology extension deformations}
Suppose that $\Lif(\xi,A')^G$ is non-empty, and let $f:\xi\ra\xi'$ be a lifting over $\Spec(A)\subset\Spec(A')$
whose isomorphism class is fixed under the $G$-action. Then the $G$-action on $\xi$ extends to an action on $\xi'$
if and only if the resulting cohomology class $[\tilde{G}]\in H^2(G,I\otimes_k\Aut_{\xi_0}(\xi_{k[\epsilon]}))$ vanishes.
\end{theorem}

In the following applications, we assume that the group $G$ is finite, and write $n=\ord(G)$ for its order.

\begin{proposition}
\mylabel{order characteristic}
Suppose   $\Lif(\xi,A')$ is non-empty and that the group order $n\geq 1$ is invertible in the residue field $k$.
Then the $G$-action on $\xi$ extends to an action on $\xi'$ for some lifting $f:\xi\ra\xi'$.
\end{proposition}

\proof
The cohomology group  $H^1(G,I\otimes_kT_{\xi_0}\shF)$ is a vector space over the field $k$,
and at the same time an abelian group annihilated by   $n=\ord(G)$. Thus it must be the zero group,
and Theorem \ref{cohomology fixed deformation} ensures that there is a lifting $\xi\ra\xi'$
over $\Spec(A)\subset\Spec(A')$
whose isomorphism class is fixed under the $G$-action.
Arguing as above, the  cohomology group $H^2(G,I\otimes_k\Aut_{\xi_0}(\xi_{k[\epsilon]}))$
vanishes, and Theorem \ref{cohomology extension deformations} tells us that
we may extend the $G$-action from $\xi$ to $\xi'$.
\qed
 
\begin{proposition}
\mylabel{sylow subgroup}
Suppose   $\Lif(\xi,A')$ is non-empty and that the residue field $k$ has  characteristic   $p>0$. 
Let $P\subset G$ be a Sylow $p$-subgroup.
Then   $\Lif(\xi,A')$ has a $G$-fixed point if and only if it has a $P$-fixed point.
Moreover, for each $[\xi']\in\Lif(\xi,A')^G$,  the $G$-action on $\xi$ extends to  $\xi'$ 
if and only if the $P$-action extends.
\end{proposition}

\proof
According to \cite{Brown 1982},  Chapter III, Proposition 10.4 the restriction map
$$
H^1(G,I\otimes_k T_{\xi_0}\shF) \lra H^1(P,I\otimes_k T_{\xi_0}\shF)
$$
is injective, and the first assertion follows from Theorem \ref{cohomology fixed deformation}. If there is
a lifting $\xi\ra\xi'$ whose isomorphism class is $G$-invariant,
we again have an injective restriction map
$$
H^2(G,I\otimes_k\Aut_{\xi_0}(\xi_{k[\epsilon]})) \lra H^2(P,I\otimes_k\Aut_{\xi_0}(\xi_{k[\epsilon]})),
$$
and the second assertion follows from  Theorem \ref{cohomology extension deformations}.
\qed

\medskip
Recall that a finitely generated free $kP$-modules $V$ have trivial cohomology groups $H^i(P,V)$,
for all $i\geq 1$. We thus get:

\begin{corollary}
\mylabel{free module}
Assumptions as in the proposition.
Then   $\Lif(\xi,A')$ has a $G$-fixed point if $\Lif(\xi,A')$ is free as $kP$-module.
Moreover, for each $[\xi']\in\Lif(\xi,A')^G$,  the $G$-action on $\xi$ extends to  $\xi'$ 
if  $Aut_{\xi_0}(\xi_{k[\epsilon]})$ is free as $kP$-module.
\end{corollary}

In some sense, this seems to be the  best possible general result:
According to \cite{Brown 1982}, Chapter VI, Theorem 8.5,
for every finite $p$-group $P$ and every field $k$ of characteristic $p>0$,
the following holds for   $kP$-modules $V$:
$$
H^1(P,V)=0\quad\Longleftrightarrow\quad H^2(P,V)=0 \quad\Longleftrightarrow\quad \text{the $kP$-module $V$ is free}.
$$
If $P$ is cyclic of order $p^\nu$ and $V$ is finitely generated, then the action of any generator $\sigma\in P$
can be viewed as a direct sum $\sigma=J_{r_1}\oplus\ldots\oplus J_{r_m}$ of
Jordan matrices $J_r\in \GL_r(k)$ with eigenvalue $\lambda=1$.
In this case, the $kP$-module $V$ is free if and only if all summands have maximal size $r_i=p^\nu$.

\section{Liftings and algebraization of finite \'etale coverings}
\mylabel{Liftings etale coverings}

Let $\Lambda$ be an adic noetherian ring, with ideal of definition $\ideala\subset \Lambda$.
In other words, the ring $\Lambda$ is noetherian, and separated and complete with respect to the $\ideala$-adic topology.
For example, the ring $\Lambda$ could be a complete local noetherian ring.
Let $Y\ra\Spec(\Lambda)$ be a proper morphism, and set $Y_0=Y\otimes_\Lambda\Lambda/\ideala$.
Write $(\Sch/Y)$ for  the category of $Y$-schemes, and 
let  $(\FEt/Y)$ the full subcategory whose objects are the  $Y$-schemes whose structure morphism 
$\pi:X\ra Y$ is finite and \'etale.
The goal of this section is to establish the following:

\begin{theorem}
\mylabel{lifting etale}
In the above situation, the restriction functor 
$$
(\FEt/Y)\lra (\FEt/Y_0),\quad X\longmapsto X\times_YY_0
$$
is an equivalence of categories.
\end{theorem}

\proof
The main task is to show that the restriction functor   is essentially surjective.
To do so, let $\pi_0:X_0\ra Y_0$ be a finite \'etale morphism.
Set $\Lambda_n=\Lambda/\ideala^{n+1}$,
and consider the infinitesimal neighborhoods $Y_n=Y\otimes_\Lambda\Lambda_n$.
According to \cite{EGA IVd}, Theorem 18.1.2, the restriction functors $(\FEt/Y_m)\ra(\FEt/Y_n)$
are equivalences of categories for all $m\geq n$.
Inductively, we choose a finite \'etale $\pi_n:X_n\ra Y_n$ and  cartesian diagrams
\begin{equation}
\label{cartesian square}
\begin{CD}
X_n		@>>>	 X_{n+1}\\
@V\pi_n VV		@VV\pi_{n+1}V\\
Y_n		@>>>	Y_{n+1}.
\end{CD}
\end{equation}
This gives a direct system $(X_n)_{n\geq 0}$ of $Y$-schemes, and in turn
a locally ringed space $\foX=(X_0,\O_{\foX})$  whose structure sheaf is
$\O_\foX=\invlim \O_{X_n}$. Then $\foX$ is a formal scheme, according to \cite{EGA I}, Proposition 10.6.3.
Let $\shI_{ji}$ be the kernels of the canonical surjections
$u_{ji}:\O_{X_i}\ra\O_{X_j}$. The closed subscheme $Y_j\subset Y_i$ corresponds to the
coherent ideal $\ideala^{i+1}\O_{Y_j}$. Since the diagrams \eqref{cartesian square} are cartesian,
we have $\shI_{ji}=\ideala^{j+1}\O_{Y_i}$.
Setting $\shI_i=\shI_{0i}$, we get $\shI_{ji}=\shI_i^{j+1}$. Thus \cite{EGA I}, Corollary 10.6.4
applies, and we infer that the formal scheme $\foX$ is adic and noetherian.

Likewise, we define $\foY=(Y_0,\O_\foY)$ via $\O_\foY=\invlim \O_{Y_n}$, which
is also an adic noetherian scheme, in fact the formal completion of $Y$ along $Y_0\subset Y$.
The diagrams \eqref{cartesian square} yield a morphism of locally ringed spaces $\pi_\infty:\foX\ra\foY$,
by \cite{EGA I}, Corollary 10.6.11.
According to \cite{EGA IIIa}, Proposition 4.8.1 this morphism $\pi_\infty:\foX\ra\foY$ is 
finite. By assumption, the formal scheme $\foY$ is algebraizable. According to \cite{EGA IIIa}, Proposition
5.4.4 the same holds for $\foX$.
In other words, there is a proper morphism $X\ra\Spec(\Lambda)$ of schemes so that our formal scheme
$\foX$ is isomorphic to the formal completion along $X_0=X\otimes_\Lambda\Lambda_0$.

The morphism $\pi_\infty:\foX\ra\foY$ comes from a unique $\Lambda$-morphism $\pi:X\ra Y$,
according to \cite{EGA IIIa}, Theorem 5.4.1.
In fact, the algebraization is the relative spectrum $X=\Spec(\shA)$ of some finite $\O_Y$-algebra $\shA$
whose formal completion becomes the finite $\O_\foY$-algebra $\shB=(\pi_\infty)_*(\O_\foX)=\invlim\O_{X_n}$.
This ensures that the morphism $\pi:X\ra Y$ is finite.
Moreover, the base-change $X\otimes_\Lambda\Lambda_n$ is isomorphic to $X_n$.
In particular, we recover the original finite \'etale covering $\pi_0:X_0\ra Y_0$.

We still have to check that the morphism $\pi:X\ra Y$ is \'etale.
By the Local Flatness Criterion (\cite{Matsumura 1980}, Section 20, Theorem 49),
the finite morphism of locally ringed spaces $\pi_\infty:\foX\ra \foY$ is flat.
It follows that the $\O_\foY$-module $\shB$ is locally free of finite rank.
In light of \cite{EGA I}, Corollary 10.8.15 the same holds for the coherent
$\O_Y$-module $\shA$. Thus $\pi:X\ra Y$ is flat.
Consider the    map
\begin{equation}
\label{associated to trace}
\shA\lra\shA^\vee,\quad a\longmapsto (a'\longmapsto\operatorname{Tr}_{\shA/\O_Y}(x\mapsto aa'x)) 
\end{equation}
associated to the trace map. The finite flat $\O_Y$-algebra $\shA$ is \'etale if and only if the trace map is surjective
(\cite{EGA IVd}, Proposition 18.2.3).
Clearly, the induced map $\shB\ra\shB^\vee$ on formal completion is associated to the trace map for $\shB=\invlim\O_{X_n}$,
and the latter is surjective. According to \cite{EGA IIIa}, Corollary 5.1.3  
the map \eqref{associated to trace} is already surjective.

Summing up, we have shown that $\pi:X\ra Y$ is finite and \'etale, and it induces the given $\pi_0:X_0\ra Y_0$.
It remains to show that the restriction functor $X\mapsto X_0$ is fully faithful.
But this follows easily from \cite{SGA 1}, Expos\'e 1, Corollary 8.4 together with
\cite{EGA IIIa}, Theorem 5.4.1.
\qed

\medskip
Note that the above result strengthens
\cite{SGA 1}, Expos\'e 1, Corollary 8.4, because it says that the finite \'etale covering
$\pi_0:X_0\ra Y_0$ not only admits a formal lifting, but even an algebraic lifting.

Now suppose that $\Lambda$ is a complete local noetherian ring with residue field $k=\Lambda/\maxid_\Lambda$,
and let $X_0$ be a proper $k$-scheme endowed with a free action of a finite group $G$.
We see that if the quotient $Y_0=X_0/G$ admits a lifting $Y\ra\Spec(\Lambda)$,
then also $X_0$ admits a lifting $X\ra\Spec(\Lambda)$, and the $G$-actions extend to
all infinitesimal neighborhoods $X_n$ and thus to $X$.
Consequently, the first  obstruction in Theorem \ref{cohomology fixed deformation} 
vanishes, and one always may choose   extensions so that the 
second obstruction in Theorem \ref{cohomology extension deformations} vanishes as well.

\section{Pullbacks of pre-Tango structures}
\mylabel{}
\def\OO{\mathcal{O}}
\let\tensor=\otimes

In this section, we observe that pre-Tango structures are preserved under finite
generically \'etale morphisms, and state some consequences 
concerning Kodaira  vanishing.

Let $k$ be an algebraically closed ground field of characteristic $p>0$.
Suppose  $Y$ is an integral smooth  projective scheme, with generic point $\eta\in Y$
and function field $F=k(Y)=\kappa(\eta)$. 
Let $\iota:\Spec(F)\ra Y$ be the inclusion of the generic point.
For arbitrary divisors $D\in\Div(Y)$, one gets  an injection
of quasicoherent sheaves $\Omega^1_{Y/k}(-D)\subset \iota_*\iota^*\Omega^1_{Y/k}$,
and thus an inclusion 
$$
H^0(Y,\Omega^1_{Y/k}(-D))\subset H^0(Y,\iota_*\iota^*\Omega^1_{Y/k})= \Omega^1_{Y/k,\eta}.
$$
An ample divisor $D\in\Div(Y)$
is called a \emph{pre-Tango} structure if there is 
a rational function $r\in F$ that is not a $p$-th power such that 
the rational differential $dr\in\Omega^1_{Y/k,\eta}$ comes from  a global section $dr\in H^0(Y,\Omega^1_{Y/k}(-pD))$.
As a short hand, one then    writes $(dr)\geq pD$.
For curves, this notion goes back to Tango \cite{Tango 1972}.
It was used by Raynaud \cite{Raynaud 1978} to construct counterexamples for Kodaira
Vanishing for surfaces, and was studied in higher dimensions by Mukai \cite{Mukai 2013}
and Takeda \cite{Takeda 2007}.

The invertible sheaf $\shL=\O_Y(D)$ associated to a pre-Tango structure $D\in\Div(Y)$ has $H^1(Y,\shL^{\otimes -1})\neq 0$,
hence is a counterexample to Kodaira Vanishing.
Conversely, if $\shL$ is ample with $H^1(Y,\shL^{\otimes -1})\neq 0$, then
$\shL^{\otimes n}$, for some integer $n\geq 1$,  comes from a pre-Tango structure $D\in\Div(Y)$,
according to \cite{Takayama 2014}, Proposition 8.

\begin{theorem}
\mylabel{pre-tango pullback}
Assume that  $X$ is another  integral smooth projective scheme and let
$\pi:X\ra Y$ be a finite surjective morphism that generically \'etale.
If $D\in\Div(Y)$ is a pre-Tango structure, then $\pi^*(D)\in\Div(X)$
is a pre-Tango structure as well.
\end{theorem}

\proof
Since $\pi:X\ra Y$ is finite, the preimage $\pi^*(D)$ of the  ample divisor $D$ remains ample.
Since this morphism is also surjective, the two schemes $X$ and $Y$ have the same dimension $d\geq 0$.
The smoothness of $X$ and $Y$ ensures that the coherent sheaves $\Omega^1_{X/k}$ and $\Omega^1_{Y/k}$ are locally
free of rank $d$.
Consider the exact sequence
$$
\pi^*(\Omega_{Y/k}^1) \lra \Omega^1_{X/k} \lra \Omega^1_{X/Y}\lra 0.
$$
The term on the right vanishes at the generic point of $X$, because the extension of function fields $k(Y)\subset k(X)$ is separable.
Consequently, the kernel $\shF$ for the map  $\pi^*(\Omega_{Y/k}^1)\ra\Omega^1_{X/k}$  on the left has rank zero.
Since the scheme $X$ has no embedded components, the kernel $\shF$ vanishes, and we get an injection
$\pi^*(\Omega_{Y/k}^1)\subset\Omega^1_{X/k}$. 
We may regard this as in inclusion inside the quasicoherent sheaf $\iota_*\iota^*\Omega_{X/k}^1$ of rational 
differentials, where $\iota:\Spec k(X)\ra X$ is the inclusion of the generic point. This gives  inclusions
$$
\pi^*(\Omega_{Y/k}^1(-pD))\subset\Omega^1_{X/k}(-p\pi^*(D))\subset \iota_*\iota^*\Omega_{X/k}^1.
$$
Since $D\in\Div(Y)$ is a pre-Tango structure, there is a rational function $r\in k(Y)$ that is not a $p$-th power
with $(dr)\geq pD$. Using that $k(Y)\subset k(X)$ is separable, one easily infers that $r\in k(X)$ is not a $p$-th power.
The rational differential $dr\in\Omega^1_{Y/k,\eta}$ extends to a global section $dr\in H^0(Y,\Omega^1_{Y/k}(-pD))$.
Viewed as element in $H^0(X,  \iota_*\iota^*\Omega_{X/k}^1)$, it lies in 
$$
H^0(X, \pi^*(\Omega^1_{Y/k}(-pD)))\subset H^0(X,\Omega_{X/k}^1(-p\pi^*(D))\subset H^0(X,  \iota_*\iota^*\Omega_{X/k}^1).
$$
Therefore, $\pi^*(D)\in\Div(X)$ is a pre-Tango structure.
\qed

\medskip
Let us say that \emph{$H^1$-Kodaira vanishing} holds on $X$ if $H^1(X,\shL^{\otimes -1})=0$ 
for all ample invertible sheaves on $\shL$. In other words, there are no pre-Tango structures on $X$.
We record the following consequence:

\begin{corollary}
Assumptions as in the theorem.
If $H^1$-Kodaira vanishing holds for $X$, then it also holds for $Y$.
\end{corollary}

Finally, we relate our conditions to liftings over the truncated ring of Witt vectors $W_2(k)$:

\begin{corollary}
Assume that  $X$ is  an  integral smooth projective scheme of dimension $d\geq 2$, and 
$\pi:X\ra Y$ be a finite surjective morphism that generically \'etale, with $Y$ smooth.
Suppose that   that $X$ admits a lifting over $W_2(k)$.
Then  $H^1$-Kodaira vanishing holds on $Y$, and furthermore $Y$ is projective.
\end{corollary}

\proof
According to \cite{EGA II}, Proposition 6.6.1, the scheme $Y$ is projective.
The assumption $d\geq 2$ and the liftability of $X$ ensures that
$H^1$-Kodaira vanishing holds for $X$, according to Deligne and Illusie 
(\cite{Deligne; Illusie 1987}, Corollary 2.8). By the previous Corollary, $H^1$-Kodaira vanishing holds on $Y$ as well.
\qed


\end{document}